\magnification 1200

%
% My TeX style file
% 20 May 2005
%
% minor improvements 22 Dec 2005
% finite fields 4 Sept 2008
% more fonts 3 Dec 2009
% minor improvements 31 October 2013

% debugging hack to make sure file is being read
%\centerline{\bf standard.tex}
%

%\magnification 1200
\settabs 10\columns
\parskip=\medskipamount\parindent=0pt

%%%%%%%%%%%%%%%% formatting stuff %%%%%%%%%%%%%

\def\H{\medskip\hrule\medskip}

%
% for some reason when this file is input by LectureNoteFormat.tex, the
% \CB macro didn't work, so I reverted to spelling it out
%
\def\Sect#1{\H\noindent\centerline{\bf{#1}}} % obscure bug fix 6 Sept 2012
\def\subSect#1{\medskip\noindent{\bf#1.}}
     % for the majiscule impaired

    % stands for ``give 'em''

%%%%%%%%%%%%%%%%%%% useful hacks %%%%%%%%%%%%

%%%%%%%%%%%%%%%%%%% math stuff %%%%%%%%%%%
\def\Cn#1{{\bf C}^{#1}}

\def\intprob#1,#2,#3{$ {\displaystyle \int_{#1}^{#2} #3 }$}

%%%%%%%%%%%%%%%%% font stuff %%%%%%%%%%%%%%
\font\small=cmr7

   % [for arXiv] /Users/wolpjame/TeX/standard.tex
\input epsf

%%%%%%%%%%%%%%%%%%%%%%

\medskip

\centerline{\bf Period Constraints on Hyperelliptic Branch Points}
\centerline{James S.~Wolper}
\centerline{Department of Mathematics}
\centerline{Idaho State University}
\centerline{Pocatello, ID~~83209--8085~~USA}
\centerline{29 June 2021}
% 2020

\medskip
\parindent=32pt
\indent
\vbox{\small\hsize=300pt
\noindent{\bf Abstract.} 
Information Theoretic analysis of the periods of a hyperelliptic curve
provides more information about the well--known but abstract
relationship between the branch points and the periods.  Here one constructs a canonical
homology basis for a hyperelliptic curve that shows that its
periods must satisfy certain constraints and defines an open set in the
Siegel upper half space that cannot contain any period matrices
of hyperelliptic  curves.
}

\parindent=0pt

%%%%%%%%%%%%%%%%%%%%%%%%%%%%%%%%%%%%%
\Sect{Introduction}

For decades now researchers have used algebraic curves to address questions 
in coding and cryptography.  This
raises interesting questions about the curves themselves based on Information
Theoretic considerations.  In one direction, 
suppose that Alice wants to tell Bob about a compact curve
({\it ie\/}, Riemann Surface) with genus $g \ge 2$.
By Torelli's Theorem, she need only send him the period matrix, 
transmitting $O(g^2)$ complex numbers.
The curve she is describing only depends on $O(g)$  parameters,
the dimension of the moduli space being $3g-3$,
so there is a lot of redundancy in her message.
The period matrix is {\sl sparse\/} in the sense of [D], and should therefore be compressible.  

The perspective that the period matrix is a compressible signal is
the central idea of the {\sl Information--Theoretic Schottky Problem\/}.
The attempt to apply ideas from Information Theory and  Compressed Sensing [D] to the
Schottky problem has led to many interesting experiments, conjectures,
and theorems [W12].  

These questions only make sense when the genus is large.
While considering them don't think of a compact surface of
genus 3 or 4 or so; the more likely model is below:

\epsfysize=3.05in
\centerline{\epsfbox{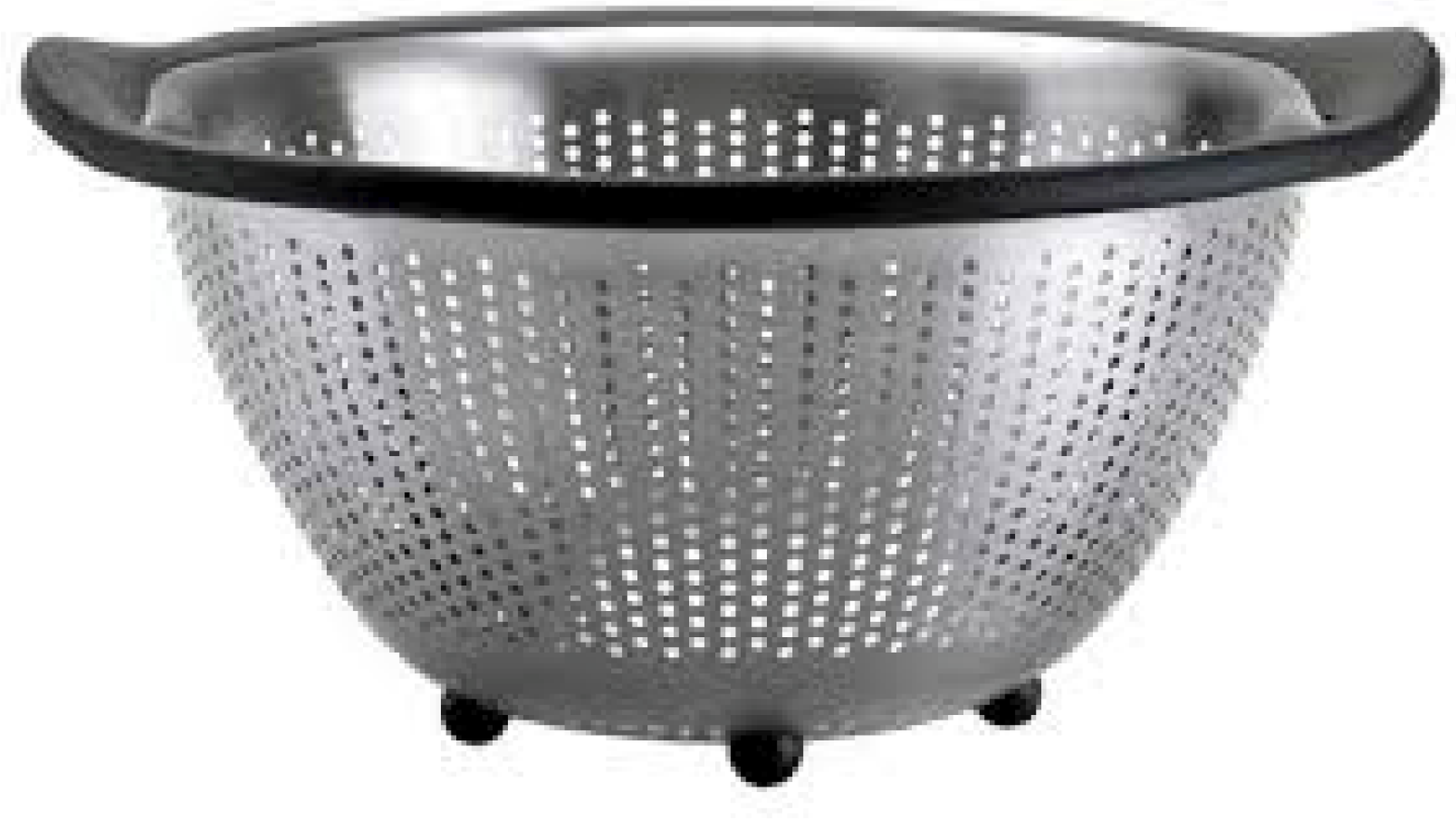}}

\centerline{{\sl Fig.~1: Compact Riemann Surface of Large Genus\/}}
\smallskip

The result described here is purely analytical, rather than 
computational; however, it was inspired by an attempt to implement ideas
in blind Compressed Sensing, as described in [GE].  

%%%%%%%%%%%%%%%%%%%%%%%%%
\subSect{Period Matrices}

For an introduction to periods see [GH].

The standard constructions for the period matrix of a compact
Riemann Surface (henceforth {\sl curve\/})
$X$ of positive genus $g$ go as follows.  Choose 
a symplectic basis $\alpha_1, \ldots, \alpha_g$, $\beta_ 1,
\ldots, \beta_g$ for the singular homology $H_1(X, {\bf Z})$;
this means that the intersections
$\alpha_i\cdot\alpha_j$,
$\beta_i\cdot \beta_j$, and $\alpha_i \cdot \beta_j$
are zero when $i \ne j$, and $\alpha_i\cdot \beta_j = 1$.
Then choose 
a basis $\omega_1,
\ldots, \omega_g$ for the space $H^{1,0}(X)$ of holomorphic differentials
on $X$, normalized
so $\int_{\alpha_i} \omega_j = \delta_{ij}$, the Dirac delta.  The
matrix $\Omega_{ij} := \int_{\beta_i} \omega_j$ is the {\sl period
matrix\/}; Riemann proved that it is symmetric with positive definite imaginary
part.  The torus $\Cn g / [I | \Omega]$ is the {\sl Jacobian\/} of $X$ ($I$ is
the $n\times n$ identity matrix).
By Torelli's Theorem, the Jacobian determines all of the
properties of $X$. In practice deciding which properties apply
is seldom successful (but see [W07]).

The normalized period matrix, whose left half is
the $g\times g$ identity, is symmetric with positive-definite imaginary part,
and the space of such matrices forms the {\sl Siegel upper half-space\/} ${\cal H}_g$.
Its dimension is $g(g+1)/2$, while the moduli space of 
curves of degree $g$ has dimension $3g-3$.  Distinguishing the period
matrices from arbitrary elements of ${\cal H}_g$ is the {\sl Schottky
Problem\/}.  See [G] for details on the problem and some of its
previous solutions.

Information--Theoretic or statistical analysis of the periods
depends on constructing a set of real numbers from the periods; typical
choices are magnitude--squared and argument.  The primary tool is then to 
sort the  modified periods in descending or, in the
extreme, non--decreasing order.  It is possible to construct

For example, the magnitudes of the periods of the Fermat curve of 
degree 9 have the distribution in Figure 1.  This data was generated using
Maple's {\tt PeriodMatrix} routine [De], then sorted using
Microsoft Excel.  The process is further explained below.

\smallskip
\epsfysize=3in
\centerline{\epsfbox{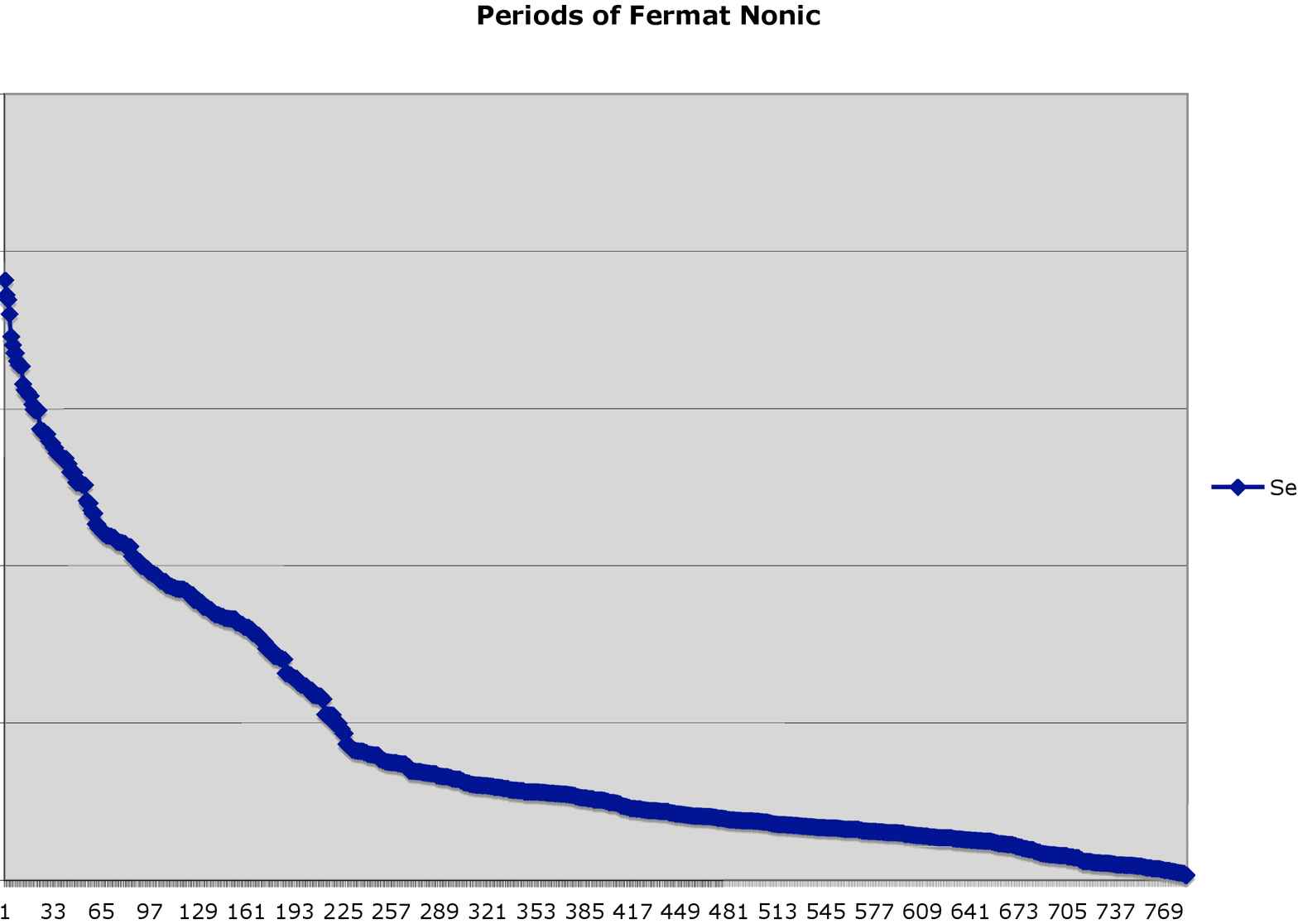}}
\centerline{{\sl Fig.~2: Periods of the Fermat Curve\/}}
\smallskip

It is possible to construct an Abelian variety with any non-decreasing distribution
by putting the periods into a symmetric matrix and adjusting so the 
imaginary part is positive--definite.  In numerical
experiments distributions from curves always have the concavity 
suggested in Figure 1.

%%%%%%%%%%%%%

\subSect{Hyperelliptic Curves}

Hyperelliptic Curves admit a degree 2 cover of the Riemann Sphere;
if the genus is $g$ one has an equation
$$
\displaystyle y^2 = \prod_{i=0}^g (x - P_i)(x - Q_i) = f(x)
$$
the curve ramifies at the $P$'s and $Q$'s.  (These are assumed to be
distinct.)

It is well--known [FK] that the branch points are holomorphic
functions of the periods; this is an abstract result and there is no
information about how the distribution of the periods
affects the placement of the branch points.  See [R] for 
an indication of the difficulties involved.  Nor are there many results
on finding properties of a curve from its period
matrix, which Torelli's Theorem suggests is possible; [W07] is one example.

The main result is:

\subSect{Theorem} It is not possible to construct a hyperelliptic curve with an
arbitrary period distribution.

The proof is by showing that it is impossible for a hyperellliptic curve
to have all of its periods equal.  This is a fairly simple argument
based on careful choice of a canonical homology basis and examination of
the integrals that define the periods.

Note that this result is consistent with earlier numerical
results [W12] that suggest that the periods of
 hyperelliptic curve tend to have little variance in their arguments,
 that is, the periods are clustered near a line in the complex plane.

%%%%%%%%%%%%%%%%%%%%%%%%%%%%%%%%%%%%%%%%%
\Sect{Period Distributions}

The matrix structure of the periods 
is only relevant here in that it has to have a 
positive--definite imaginary part.  But because a permutation of the 
bases of $H_1(X, {\bf Z})$ or $H^{1,0}(X)$ permutes the matrix entries
they are not intrinsic so one works with {\sl lists\/} of periods.

A commonly--used strategy for signal compression is to expand
the signal in some kind of series (such as a Fourier series)
and use only the terms with the ``biggest" coefficients.  This
has also long been a strategy with Singular Value Decompositions;
see [GD].  (This is distinct from the noise--reduction
strategy which uses only the lower--frequency terms.)   It is
not widely appreciated that E.~Lorenz used this technique
to introduce Ordinary Differential Equation approximations
to a set of Partial Differential Equations in the first
papers on Chaos [L].

When thinking of periods as a distribution or signal,
then, the first step is to 
create a set of real numbers from the periods.  Various
methods have been used: the modulus, the squared--modules,
and the argument.
Now sort the real versions of the periods in 
descending order to make a list $p_1$,
$p_2, \ldots, $ 
$p_{g(g-1)/2}$, and plot $(n, p_n)$.

 If the periods were uniformly distributed
the plot would be a straight line, but numerical 
experiments suggest that it is concave up.

A symmetric matrix with positive--definite imaginary part is
the period matrix of an abelian variety: the theta--function constructed
from it converges.  
The Schottky Problem comes into play by noting that any {\sl random\/}
set of periods, perhaps adjusted by a few factors of the form $e^{i\phi}$
(where $\phi$ is a real number)
to make the imaginary part of the matrix positive--definite,
defines an abelian variety.
In particular,  it is possible to define an {\sl abelian variety\/} all of whose
periods have the same modulus, but, conjecturally,
this is not the Jacobian of a Riemann Surface.

%%%%%%%%%%%%%%%%%%%
\Sect{The Canonical Basis}

The analysis of the periods uses a specific canonical homology basis,
which is best explained pictorially for genus 3.  The branch points are labelled
$P_i$, $Q_i$ for $i = 1, 2, 3$.  It is convenient
to have one branch point at $\infty$.

\epsfxsize=2.05in
\centerline{\epsfbox{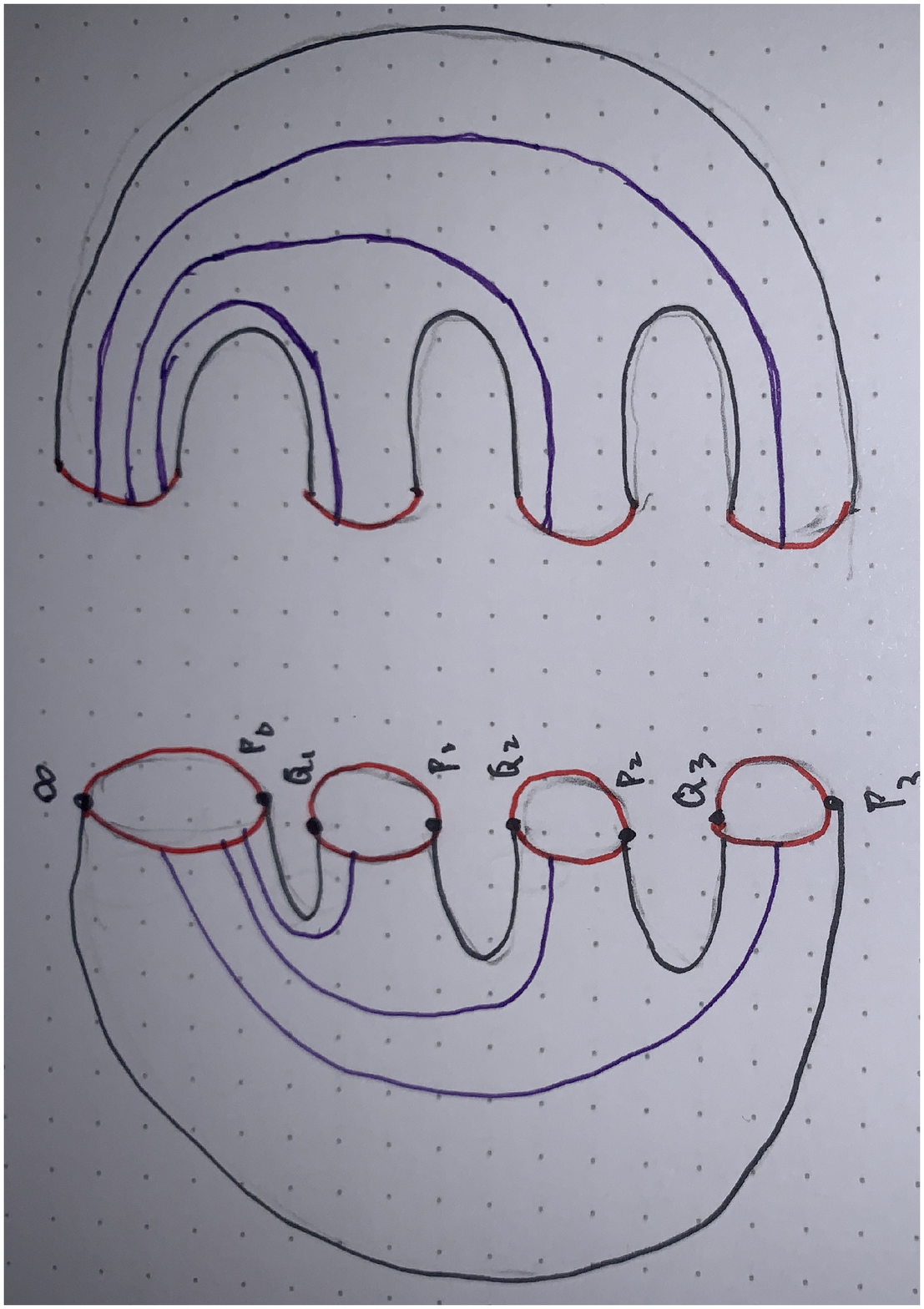}}

\centerline{{\sl Fig.~3: Canonical Homology Basis\/}}
\smallskip

Figure 3 omits the orientations of the cycles for clarity.  The
double cover of the Riemann sphere is evident, and the
``tubes" represent the effect of crossing a cut conceptually,
not literally.

Figure 4 details the cuts in the target plane of the 
hyperelliptic projection.  Dashed portions of the cycles are on opposite sheets,
so an evident intersection between a dashed cycle and a solid cycle
is not an actual intersection

\epsfysize=3.05in
\centerline{\epsfbox{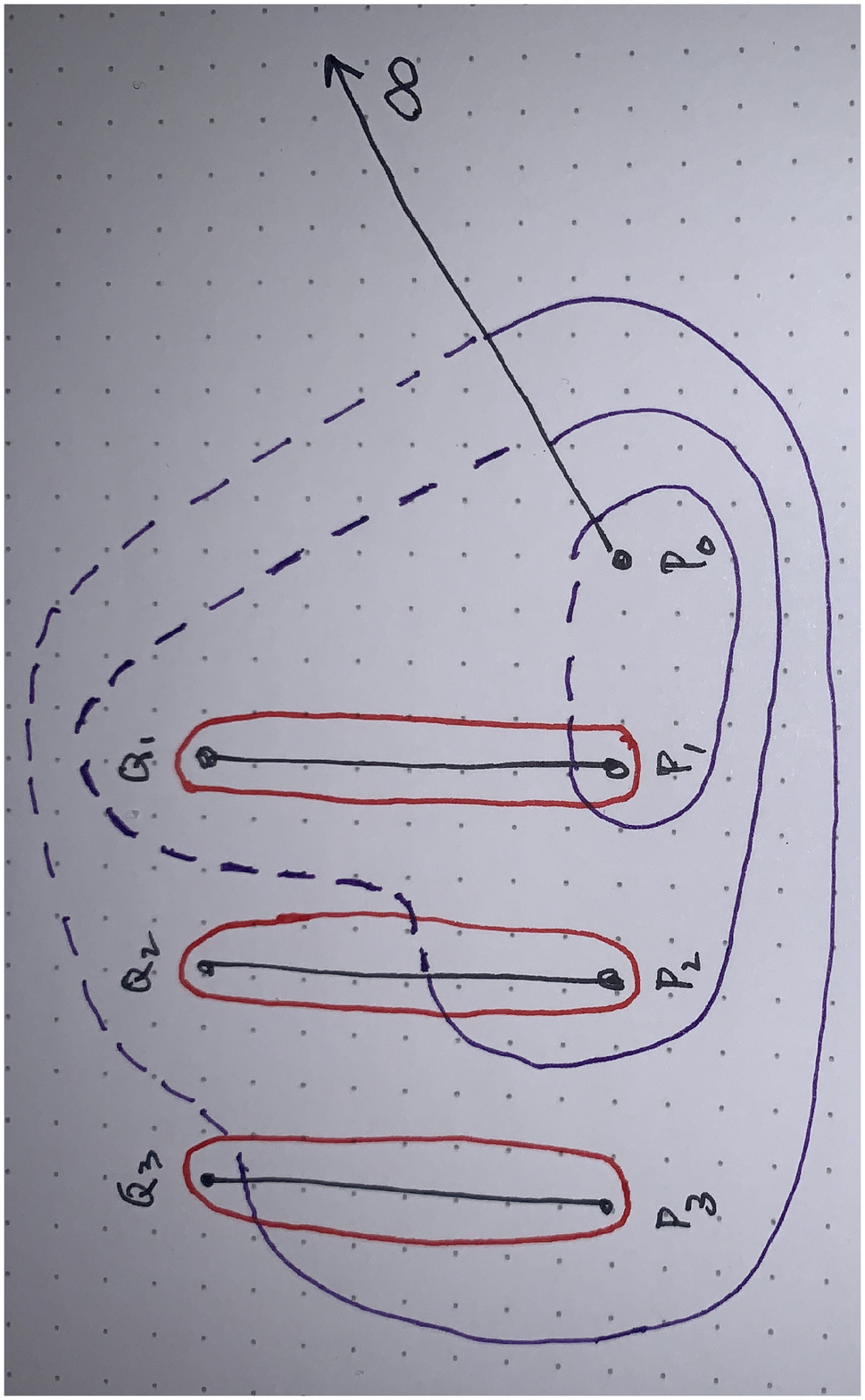}}

\centerline{{\sl Fig.~4: Canonical Homology Basis\/}}
\smallskip

%%%%%%%%%%%%%%%%%%%%%%%%%
\Sect{Branch Point Constraints}

To determine the periods one integrates the cohomology basis over
each of the cycles in the canonical homology basis.  Refer to
Figure 4: there are three cycles, each surrounding the pairs $P_i$, $Q_i$.
Call these $\alpha_i$.  

The rest of the argument is quite simple: in the picture, there is a
homology between $\alpha_1 + \alpha_2$ and $\alpha_3$,
so for any  holomorphic 1--form $\omega$
$$
\int_{\alpha_1}\omega + \int_{\alpha_2} \omega =
\int_{\alpha_3} \omega;
$$
similar statements hold in general for genus larger than 3.
But if the periods are equal then, in this case, one of the
integrals is zero, which is absurd.

%%%%%%%%%%%%%%%%%%%%
\Sect{Conclusion}

The significance of this result is in the consequence, not the
proof of the theorem.  It identifies a locus in
the moduli space of hyperelliptic
surfaces which are not jacobians, thus determining a 
partial solution to the Schottky problem.  Also
significant is the information--theoretic 
perspective that led to this result.

%%%%%%%%%%%%%%%%%%%%%%%%%%%%%%%%%%%%%%%%%
\Sect{References} 

\medskip \frenchspacing

\noindent

\item{[De]} Deconinck, B.~and Mark van Hoeij, ``Computing Riemann
Matrices", {\sl Physica D\/} {\bf 152-153} (2001), 28 -- 46.

\item{[D]} Donoho, D.,  ``Compressed Sensing," {\sl IEEE Transactions
on Information Theory\/}, {\bf 52} (April, 2006), no.~3, 1289 -- 1306.

\item{[FK]} Farkas, H., and I.~Kra, {\sl Riemann Surfaces\/}.  NY: Springer (1991).

\item{[G]} Grushevsky, Samuel, ``The Schottky Problem," 
{\tt arXiv}:1009.0369v1, 2 September 2010.

\item{[GD]} Gavish, Matan and David L.~Donoho,
{\sl The Optimal Hard Threshold for Singular Values is $4/\sqrt(3)$\/},
{\tt arXiv}:1305.5870v3, 24 May 2013.

\item{[GE]} Gleicham, S., and Y.~Eldar,
``Blind Compressed Sensing," {\tt arXiv}:1002.2586v2,
28 April 2010.

\item{[GH]} Griffiths, P., and J.~Harris,
{\sl Principles of Algebraic Geometry\/}.
NY: Wiley (1979).

\item{[L]} Lorenz, E., ``Deterministic Nonperiodic
Flow," {\sl J.~Atmos.~Sci.\/} {\bf 20} March, 1963.

\item{[R]} Rauch, H., ``On the Transcendental Moduli of
Algebraic Riemann Surfaces," {\sl Proceedings of
the National Academy of Sciences\/}, {\bf 41} (1955), 42 -- 9.

\item{[W07]} Wolper, J.,  ``Analytic Computation of some
Automorphism groups of Riemann Surfaces," 
{\sl Kodai Mathematical Journal\/}, {\bf 30} (2007), 394 -- 408.

\item{[W12]} Wolper, J., ``Numerical experiments related to the \hbox{ information-theoretic }
Schottky and Torelli problems," {\sl Albanian Journal of Mathematics\/}
[Special issue on computational algebraic geometry], {\bf 6}, No.~1, (2012).

%%%%%%%%%%%%%%%%%%%%%%%%%%%%%%%%%%%%%
\bye